\documentclass[preprint,12pt]{elsarticle}

\usepackage{graphicx}

\usepackage{amssymb}
\usepackage{amsmath}
\usepackage{mathrsfs}
\usepackage{stmaryrd}
\usepackage{dsfont}

\newcommand{\ie}{\textit{i.e.} }
\newcommand{\eg}{\textit{e.g.}}

\graphicspath{{./pics/}}
\journal{ }

\begin{document}

\begin{frontmatter}



\title{Scalable computation of intracellular metabolite concentrations}


\author[label1]{Amir Akbari}
\ead{amakbari@ucsd.edu}
\author[label1,label2]{Bernhard O. Palsson}
\ead{palsson@ucsd.edu}

\address[label1]{Department of Bioengineering, University of California San Diego, La Jolla, CA 92093}
\address[label2]{Novo Nordisk Foundation Center for Biosustainability, Technical University of Denmark, 2800 Lyngby, Denmark}

\begin{abstract}
Current mathematical frameworks for predicting the flux state and macromolecular composition of the cell do not rely on thermodynamic constraints to determine the spontaneous direction of reactions. These predictions may be biologically infeasible as a result. Imposing thermodynamic constraints requires accurate estimations of intracellular metabolite concentrations. These concentrations are constrained within physiologically possible ranges to enable an organism to grow in extreme conditions and adapt to its environment. Here, we introduce tractable computational techniques to characterize intracellular metabolite concentrations within a constraint-based modeling framework. This model provides a feasible concentration set, which can generally be nonconvex and disconnected. We examine three approaches based on polynomial optimization, random sampling, and global optimization. We leverage the sparsity and algebraic structure of the underlying biophysical models to enhance the computational efficiency of these techniques. We then compare their performance in two case studies, showing that the global-optimization formulation exhibits more desirable scaling properties than the random-sampling and polynomial-optimization formulation, and, thus, is a promising candidate for handling large-scale metabolic networks.    
\end{abstract}


\end{frontmatter}






\section{Introduction}\label{sec:introduction}
Synthetic biology and metabolic engineering have significantly advanced in the past decade, creating enormous opportunities to tackle challenging problems in manufacturing and medicine \cite{Khalil2010, xie2018designing}. Tools and techniques have been developed to build molecular modules, providing control over gene-expression regulations (at the transcription, translation, and post-translation level), protein function, and inter-cellular communications \cite{Purnick2009}. These methods have been implemented in practice to design, construct, and optimize synthetic pathways to produce novel and/or high-value products \cite{lee2012systems}. Application areas abound, ranging from therapeutics and synthetic ecosystems to biofuels and biomaterials \cite{Khalil2010}. These developments led to optimal designs for heterologous pathways to synthesize 1,4-butanediol from renewable substrates in \emph{Escherichia coli} \cite{Yim2011}, synthetic pathways for the antimalarial-drug precursor isoprenoid in \emph{E. coli} and yeast \cite{martin2003engineering, ro2006production}, engineered quorum-sensing
systems for population control in microbial communities \cite{you2004programmed}, and synthetic regulatory circuits to enhance productivity and toxin tolerance in oleaginous strains for biofuel production \cite{Berla2013}. Such enabling technologies could significantly impact the future of pharmaceutics, agriculture, and energy systems, furnishing economically efficient and environmentally friendly routs to manufacture commodity chemicals.       

Synthetic-pathway design based on intuitive principles and laboratory tools is applicable to small modules with limited functionality. Therefore, predictive models and computational tools are needed to systematically evaluate alternative scenarios and the outcomes of optimization strategies at the system level \cite{Purnick2009}. In this regard, constrained-based models of metabolism (M-models) have proven useful. OptKnock \cite{Burgard2003} is an early example of model-driven strain-design algorithm that identifies optimal reaction-deletion strategies for the growth-coupled production of target metabolites using flux-balance analysis (FBA). OptORF \cite{Kim2010} is an extension of this algorithm, providing gene-deletion and overexpression strategies that lead to the overproduction of desired biochemicals using an FBA of integrated regulatory and metabolic networks. Recently, capabilities of metabolism-expression models (ME-models) in predicting byproduct secretion were leveraged to establish robust strain designs among numerous scenarios furnished by metabolism-based algorithms \cite{Dinh2018}. The predicted performance of optimized knockout mutants have been validated experimentally in several cases \cite{Dinh2018}. However, the foregoing algorithms may overlook part of regulatory responses to gene deletions because their underlying mathematical models are still far from comprehensive. Specifically, M- and ME-models do not directly account for intracellular concentrations to ensure the consistency of metabolic fluxes and thermodynamic constraints. Several key reactions in central carbon metabolism are allosterically regulated by metabolite levels with important implications for pathway design and optimization \cite{Ellis2012}.                   

This paper concerns a constrained-based treatment of intracellular concentrations. Metabolite levels inside the cell are subject to fundamental and evolutionary constrains, such as those imposed by thermodynamic laws, electroneutrality, osmotic balance, and ion homeostasis. The constrain-based framework we introduce seeks to (i) identify the conditions, under which the biophysically concentration solution space (CSS), defined by these constraints, is non-empty, and (ii) derive characteristics of the CSS, from which biological insights can be drawn. Thermodynamic constraints are represented by nonlinear polynomial inequalities with respect to concentrations; all the other constraints are linear equalities. This CSS is generally nonconvex and disconnected. Therefore, computationally tractable algorithms are necessary to study CSSs for metabolic networks of biological interest. Addressing these goals poses the same challenges as those of feasibility and optimality identification in nonlinear programs (NLPs) \cite{nocedal2006numerical}.

Determining the feasibility of NLPs is a difficult task with no guaranteed algorithm for general cases \cite{Chinneck2007}. Penalty methods are the most common approaches, where feasible points are identified by minimizing an unconstrained penalty function, measuring the extent to which the constraints of the original problem are violated \cite{elwakeil1995methods}. Bootstrapping is another approach used by several NLP solvers \cite{debrosse1973feasible, Chinneck2007}. Here, the search for feasibility starts at an initial point, where a subset of constraints are satisfied. The initial point is iteratively adjusted to reduce the number of violated constraints without affecting the set of satisfied constraints. How well these methods can converge to a feasible point heavily depends on the problem structure and quality of the initial point. Therefore, they are usually complemented by heuristic techniques for initial-point placement and low-cost preprocessing solvers for feasibility approximation \cite{Chinneck2007}.

Characteristic measures of the CSS can provide valuable information about the functional states of biological networks. The volume measure is a natural choice for characterizing the CSS (and the feasible set of NLPs in general). Although it may not be directly related to biologically interpretable parameters, it can be linked to a large class of probability measures that quantify the stochastic nature of intracellular dynamics. However, deterministic computation of the volume, exactly or approximately, is a major challenge, which cannot be accomplished in polynomial time, even for convex bodies \cite{Simonovits2003}. Nevertheless, given a failure probability, stochastic methods of polynomial complexity exist to approximate the volume of convex bodies to within a tolerable error \cite{dyer1991random}. Therefore, random-sampling techniques are commonly used to approximate the volume of convex and nonconvex bodies \cite{Simonovits2003}. Here, generating random sample points with a uniform probability distribution in irregularly shaped or nonconvex regions is a major issue \cite{Chinneck2007}. In these cases, random sampling is performed in convex enclosures using Markov Chain Monte Carlo methods, such as hit-and-run, where the probability distribution over the original set iteratively converges to uniformity at a geometric rate with a polynomial complexity in dimension \cite{Smith1984, Belisle1998}. These techniques perform poorly in higher dimensions or for long and thin regions (either the enclosure or original set) \cite{Chinneck2007}.               

Given that the constraints considered in this work are all representable as polynomials in concentrations, methods of algebraic geometry can be applied to determine the feasibility and characteristics of the CSS \cite{bochnak2013real}. A number of \emph{Positivstellens\"atze} have been proven in the literature, providing infeasibility certificates for semialgebraic sets \cite{Stengle1974, Schmuedgen1991, putinar1993positive}. All these results rely on a connection between feasibility and the existence of sum-of-squares (SOS) polynomials. Although exact, computational verification of these certificates is not straightforward. Tractable algorithms based on semidefinite-programming (SDP) relaxation of polynomial optimization problems enable the computation of infeasibility certificates in low-dimensional cases \cite{Lasserre2001, Parrilo2003}. Similar approximation algorithms have been proposed to compute the volume of basic semialgebraic sets \cite{Henrion2009}. Even though the performance of these algorithms does not rely on the convexity or connectedness of the set, their applicability is limited by their undesirable scaling properties.       

In this paper, we develop tractable computational methods to characterize the CSSs arising from a constrained-based model of intracellular metabolite concentrations, which we previously developed \cite{Akbaria2020}. We examine three approaches, including polynomial optimization, random sampling, and global optimization. We then compare the performance and scalability of these techniques through case studies. For polynomial optimization problems, we derive infeasibility certificates, which can be computationally verified by solving a hierarchy of SDPs, growing in size at higher levels. We leverage the sparsity of these SDPs, which stems from the sparse stoichiometry matrices of biochemical reaction networks, to improve the computational efficiency in lower dimensions. Finally, we introduce a penalty method, tailored to the structure of the biophysical constraints in our the model, to determine the emptiness of the CSS. This is accomplished through a phase-I formulation, which can be solved efficiently in higher dimensions using global optimization solvers. These techniques have more desirable scaling properties than polynomial-optimization approaches and can handle metabolic networks of core-scale size.      

\section{A constraint-based model of intracellular concentrations}\label{sec:model}
In this section, we introduce a simplified version the constraint-based model we recently developed \cite{Akbaria2020} to characterize the CSS for bacterial organisms. We consider four major biophysical constraints pertinent to biological systems, namely electroneutrality, osmotic balance, fixed buffer capacity, and mass-action principle. We seek to quantify the cytoplasmic concentrations for any given set of extracellular conditions. These constraints are mathematically formulated as 
\begin{equation}
\sum_{i=1}^{n}z_iC_i=0,
\label{eq:eq1}
\end{equation}
\begin{equation}
C_{t,c}-C_{t,p}+2\left[I_c f^{\phi}(I_c)-I_p f^{\phi}(I_p)\right]=\frac{\Delta\Pi}{RT},
\label{eq:eq2}
\end{equation}
\begin{equation}
\sum_{i=1}^{n}\beta_iC_i=\mathcal{B},
\label{eq:eq3}
\end{equation}
\begin{equation}
\sum_{i=1}^{n}C_i=C_{t,c},
\label{eq:eq4}
\end{equation}
\begin{equation}
\Delta_r G'^{\circ}_j+RT\ln\Gamma_j\le 0, \quad \Gamma_j:=\prod_{i=1}^{n}\left(\frac{C_i}{C^{\circ}}\right)^{S_{ij}}, \quad j=1\cdot\cdot m
\label{eq:eq5}
\end{equation}
with $C$ the vector of cytoplasmic metabolite concentrations, $C^{\circ}$ the reference concentration, $C_{t,c}$ total cytoplasmic concentration, $C_{t,p}$ total periplasmic concentration, $I_c$ cytoplasmic ionic strength, $I_p$ periplasmic ionic strength, $\mathcal{B}$ total cytoplasmic buffering capacity, $\Delta\Pi:=\Pi_c-\Pi_p$ osmotic-pressure differential, $R$ gas universal constant, $T$ temperature, $S$ stoichiometry matrix, $f^\phi$ Pitzar's function for osmotic coefficient \cite[Eq.~(5)]{pitzer1974thermodynamics}, $\Delta_r G'^{\circ}$ standard transformed Gibbs energy of reaction, $\Gamma$ reaction quotient, $z$ metabolite charge, $\beta$ metabolite buffer intensity, $j$ reaction index, $i$ metabolite index, $n$ number of metabolites, and $m$ number of reactions. We simplify Eq.~(\ref{eq:eq2}) by linearizing it around concentrations $C_{0,i}$, where the osmotic pressure differential $\Delta\Pi_0$ is known, to arrive at
\begin{equation}
\sum_{i=1}^{n}\phi_iC_i=\frac{\Delta\Pi}{RT}-C_{t,0},\quad C_{t,0}:=\frac{\Delta\Pi_0}{RT}-\sum_{i=1}^{n}\phi_iC_{0,i},
\label{eq:eq6}
\end{equation}
where 
\begin{equation}
\phi_i:=1+2\left.\frac{\partial[I_c f^{\phi}(I_c)]}{\partial C_i}\right|_0.
\label{eq:eq7}
\end{equation}
Here, the subscript $0$ indicates that the derivatives are evaluated at $C_0$---the vector of cytoplasmic concentrations, at which the osmotic pressure differential is given. An equivalent form for Eq.~(\ref{eq:eq5}) can be derived based on the strict monotonicity of the exponential function 
\begin{equation}
\prod_{i=1}^{n}\left(\frac{C_i}{C^{\circ}}\right)^{S_{ij}}\le K'_j,\quad K'_j:=\exp\left(-\frac{\Delta_rG'^{\circ}_{j}}{RT}\right), \quad j=1\cdot\cdot m,
\label{eq:eq8}
\end{equation}
where $K'_j$ is the apparent equilibrium constant of reaction $j$ \cite{Alberty2005}. This form can be readily represented as a polynomial in concentrations, so it is more suitable for polynomial optimization methods, which will be discussed in the next section. The CSS is defined as
\begin{equation}
\mathcal{C}:=\left\{C\in \mathds{R}^{n}_{+}:\mathrm{Eqs.}\; \mathrm{(\ref{eq:eq1}), (\ref{eq:eq3})-(\ref{eq:eq6})} \; \mathrm{hold}\right\}.
\label{eq:eq9}
\end{equation}

We reformulate the biophysical constraints with respect to mole fractions $x_i$, which can be represented more compactly as
\begin{equation}
\begin{cases}
Ax=w+F\theta, \\
S^{\mathrm{T}}\ln x \leq \kappa+\nu\ln\theta_1, \\
x\geq 0,
\end{cases}
\label{eq:eq10}
\end{equation} 
where $x_i:=C_i/C_{t,c}$, $\kappa_j:=\ln\left[K'_j(C^{\circ}/C_s)^{\nu_j}\right)]$, $\nu_j:=\sum_{i=1}^{n}S_{ij}$, $C_s:=\Delta\Pi/RT-C_{t,0}$, and 
\begin{equation}
A:=
\begin{bmatrix}
	z^{\mathrm{T}}\\
	\phi^{\mathrm{T}}\\
	\beta^{\mathrm{T}}\\
	1_{1\times n}
\end{bmatrix} 
,\quad w:=
\begin{bmatrix}
	0\\
	0\\
	0\\
	1
\end{bmatrix}
,\quad F:=
\begin{bmatrix}
	0 & 0\\
	1 & 0\\
	0 & 1\\
	0 & 0
\end{bmatrix}
,\quad \theta:=\frac{1}{C_{t,c}}
\begin{bmatrix}
	C_s\\
	\mathcal{B}
\end{bmatrix}
\label{eq:eq11}
\end{equation}
with $S\in\mathds{R}^{n\times m}$, $A\in\mathds{R}^{\ell\times n}$, $w\in\mathds{R}^{\ell}$, $F\in\mathds{R}^{\ell\times n_{\theta}}$, $\theta\in\mathds{R}^{n_{\theta}}$, $x\in\mathds{R}^{n}$, $\kappa\in\mathds{R}^{m}$, and $\nu\in\mathds{R}^{m}$. The thermodynamic inequalities in Eq.~(\ref{eq:eq10}) can be written in the polynomial form 
\begin{equation}
K''_j\prod_{i=1}^{n}x_i^{S_{ij}^{-}}-\prod_{i=1}^{n}x_i^{S_{ij}^{+}}\ge 0,\quad j=1\cdot\cdot m,
\label{eq:eq12}
\end{equation} 
where $K''_j:=K'_j(C^{\circ}\theta_1/C_s)^{\nu_j}$, $S_{ij}^{+}:=\max\{S_{ij},0\}$, and $S_{ij}^{-}:=-\min\{S_{ij},0\}$. The CSS in the mole-fraction space is defined accordingly 
\begin{equation}
\mathcal{X}(\theta):=\left\{x\in [0,1]^{n}:\mathrm{Eq.}\; \mathrm{(\ref{eq:eq10})} \; \mathrm{holds}\right\}.
\label{eq:eq13}
\end{equation} 
Here, the CSS is parametrically represented with respect to $\theta$ to determine part of the parameter space, for which CSS is non-empty. Whenever polynomial optimization methods are applied to study $\mathcal{X}(\theta)$, it is assumed that the thermodynamic inequalities in Eq.~(\ref{eq:eq10}) are replaced by their polynomial counterparts in Eq.~(\ref{eq:eq12}). In biological systems, intracellular concentrations can spans several orders of magnitude. Therefore, we use the logarithmic map $\Xi:=x\mapsto y=\ln x$ to reformulate Eqs.~(\ref{eq:eq10}) and (\ref{eq:eq13}) into the forms
\begin{equation}
\begin{cases}
A\exp y=w+F\theta, \\
S^{\mathrm{T}}y \leq \kappa+\nu\ln\theta_1, \\
y\leq 0,
\end{cases}
\label{eq:eq14}
\end{equation}
\begin{equation}
\mathcal{Y}(\theta):=\left\{y\in \mathds{R}_{-}^{n}:\mathrm{Eq.}\; \mathrm{(\ref{eq:eq13})} \; \mathrm{holds}\right\}.
\label{eq:eq15}
\end{equation} 
The mole-fraction and logarithmic mole-fraction spaces are referred to as the $X$ and $Y$ space, respectively. Characterizing the CSS using global optimization techniques is more computationally convenient in the $Y$ space than the $X$ space.  

\begin{figure}
\centering
\includegraphics[width=0.5\linewidth]{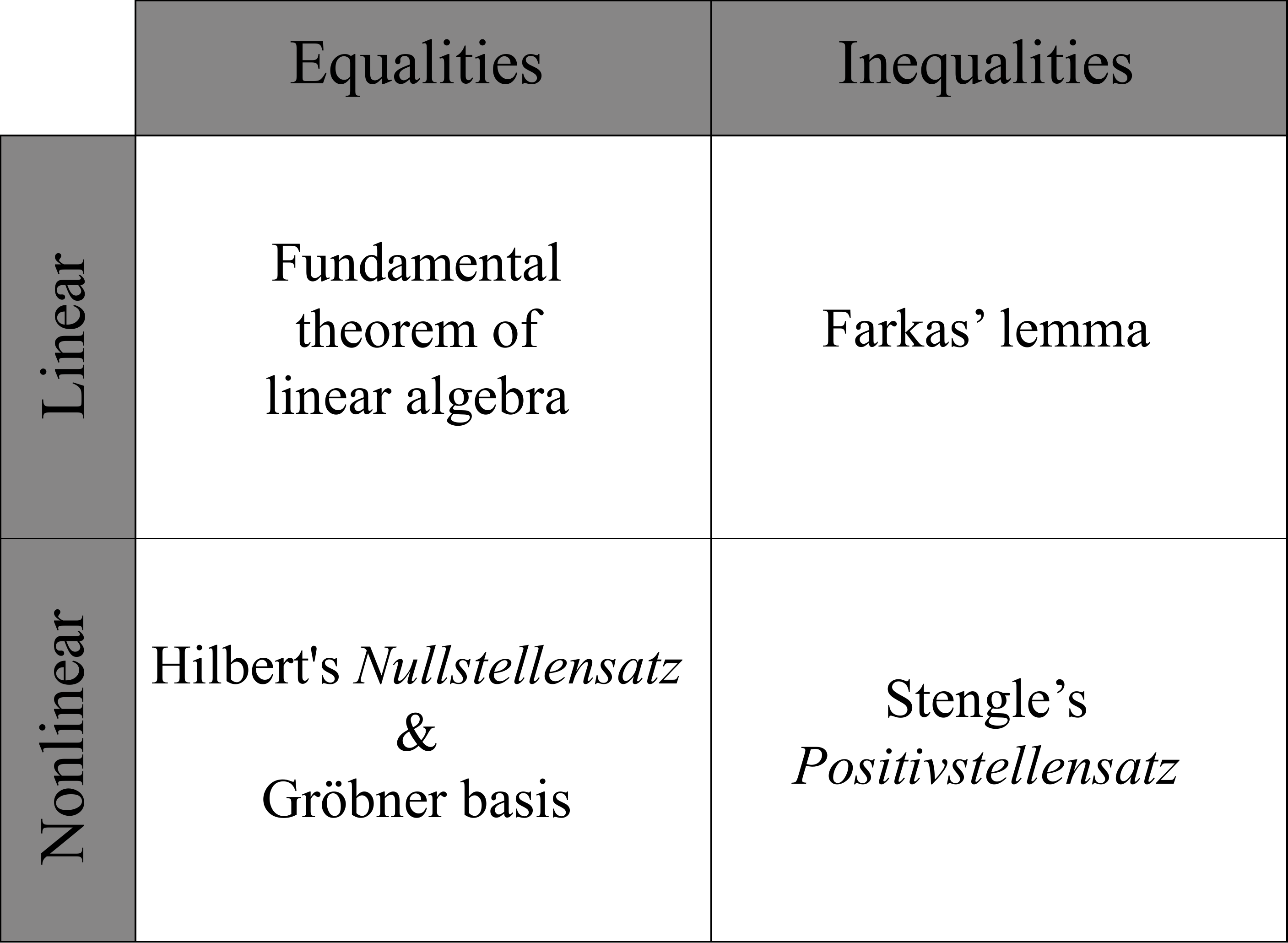}
\caption{Important theorems concerning the feasibility of semialgebraic sets classified based on their constraint type \cite{Parrilo2003,bochnak2013real}.}
\label{fig:fig1}
\end{figure}

\section{Feasibility of concentration sets}\label{sec:feasibility}
\subsection{Polynomial optimization}\label{sec:pop}
The CSS in the $X$ space is defined by a number of polynomial equalities and inequalities in $x$, and, thus, it is a semialgebraic set \cite{Parrilo2003, bochnak2013real}. Here, we consider a more general form of Eq.~(\ref{eq:eq13})
\begin{equation}
\mathcal{X}:=\left\{x\in \mathds{R}^{n}:G_j(x)\ge 0, H_i(x)=0, j=1\cdot\cdot m, i=1\cdot\cdot \ell\right\},
\label{eq:eq16}
\end{equation} 
where $H_i(x)$ are affine in $x$ and $G_j,H_i\in\mathds{R}[x]$ with $\mathds{R}[x]$ the polynomial ring over the real field. Semialgebraic sets are solution spaces of polynomial optimization problems. To determine whether an optimal solution can be found, necessary and sufficient feasibility criteria are required. Figure~\ref{fig:fig1} summarizes these criteria for several classes of constraints, by which semialgebraic sets are defined. We present a general result concerning real semialgebraic sets \cite{bochnak2013real}, which is relevant to our discussion:
\begin{equation}
\mathcal{X}=\emptyset \Leftrightarrow -1\in \mathrm{cone}\; \{G_j\}_1^m+\mathrm{span} \; \{H_i\}_1^{\ell}. 
\label{eq:eq17}
\end{equation} 
The $\mathrm{cone}$ and $\mathrm{span}$ functions in Eq.~(\ref{eq:eq17}) are defined over $\mathds{R}[x]$ similarly, but not identically, to the way they are defined over $\mathds{R}$\begin{equation}
\begin{aligned}
\mathrm{cone}\; \{G_j\}_1^m:=&\sigma(x)+\sum_{j_1=1}^{m}\sigma_{j_1}(x)G_{j_1}(x)+\sum_{j_1=1}^{m}\sum_{j_2=j_1}^{m}\sigma_{j_1j_2}(x)G_{j_1}(x)G_{j_2}(x)\\
&+\cdots+\sum_{j_1=1}^{m}\cdots\sum_{j_m=j_{m-1}}^{m}\sigma_{j_1\cdots j_m}(x)G_{j_1}(x)\cdots G_{j_m}(x),
\end{aligned}
\label{eq:eq18}
\end{equation}
\begin{equation}
\mathrm{span} \; \{H_i\}_1^{\ell}:=\sum_{i=1}^{\ell}\beta_i(x)H_i(x).
\label{eq:eq19}
\end{equation}  
Here, $\sigma,\sigma_{j_1},\sigma_{j_1j_2},\cdots,\sigma_{j_1\cdots j_m}\in\Sigma[x]$ and $\beta_i\in\mathds{R}[x]$ with $\Sigma[x]$ the SOS polynomial cone. It can be easily verified that $\mathrm{span} \; \{H_i\}_1^{\ell}$ is an ideal of $\mathds{R}[x]$ \cite{bochnak2013real}. Verification of Eq.~(\ref{eq:eq17}), which involves a search in the set of SOS and real polynomials of arbitrary degree, is an NP-hard problem and cannot be accomplished in practice \cite{Parrilo2003}. Therefore, the goal here is to find infeasibility certificates from Eq.~(\ref{eq:eq17}) through a hierarchy of optimization problems, where the search is limited to SOS and real polynomials of bounded degree. In this setting, $\Sigma[x]$ and $\mathds{R}[x]$ are approximated by $\Sigma[x]_d$ (set of SOS polynomials of degree at most $2d$) and $\mathds{R}[x]_d$ (set of real polynomials of degree at most $d$), such that $\Sigma[x]_d\to\Sigma[x]$ and $\mathds{R}[x]_d\to\mathds{R}[x]$ for arbitrarily large number of hierarchical levels. Starting from lowest degree polynomials, larger and larger subsets of $\Sigma[x]$ and $\mathds{R}[x]$ are searched for a solution to Eq.~(\ref{eq:eq17}) with every additional optimization level. At each level, SOS and real polynomials are constructed for a fixed $d$. If a solution is found, it certifies the infeasibility of $\mathcal{X}$. Because the number of thermodynamic inequalities in our model is generally larger than those typically handled by SDP relaxation methods, we further simplify our analysis by truncating the sum in Eq.~(\ref{eq:eq18}), including as many generators in the computation as necessary to find a solution at each level. 

Let $d_g$ denote the maximum total degree of all the generators of $\Sigma[x]_d$ considered in the computation, where the total degree of all the SOS polynomials is $2d$. We define 
\begin{equation}
s_p(d):={n+d\choose d}, \quad s_f(d):={n+d-1\choose d},
\label{eq:eq20}
\end{equation}  
where $s_p$ is the dimension of the set of polynomials of dimension $d$ or less. Similarly, $s_f$ is the dimension of the set of degree-$d$ forms \cite{Parrilo2003}. Let also $\rho:=2d+d_g$ and $z$ a column vector\footnote{Not to be confused with the metabolite charge in Eq.~(\ref{eq:eq1}).} containing the standard monomial basis of $\mathds{R}[x]_{\rho}$. Then, every polynomial in Eqs.~(\ref{eq:eq18}) and (\ref{eq:eq19}) can be written in matrix form \cite{Parrilo2003}
\begin{equation}
\sigma_j(x):=z^{\mathrm{T}}\mathrm{P}^jz, \quad G(x):=\mathrm{G}z,\quad H(x):=\mathrm{H}z, \quad B(x):=\mathrm{B}z.
\label{eq:eq21}
\end{equation}  
Here, $B(x):=[\beta_1(x);\cdots;\beta_{\ell}(x)]$, $\mathrm{P}^j$ are positive semidefinite matrices, and $\sigma_j(x)$ represents all the SOS polynomials in Eq.~(\ref{eq:eq18}). Moreover, $\mathrm{P}^j$, $\mathrm{G}$, $\mathrm{H}$, $\mathrm{B}$, and $z$ are matrices of dimension $s_p(\rho)\times s_p(\rho)$, $m\times s_p(\rho)$, $\ell\times s_p(\rho)$, $\ell\times s_p(\rho)$, and $s_p(\rho)\times 1$ respectively. 

To derive a matrix-form representation of Eq.~(\ref{eq:eq17}), we construct the multiplication table of $\mathds{R}[x]$. First, we define a map between the index of a monomial $u\in\mathds{R}[x]$ in the standard basis of the ring and its exponent vector $\alpha$
\begin{equation}
\mathcal{J}:\mathbb{N}^n\mapsto\mathbb{N}:\alpha\mapsto\mathcal{J}(\alpha).
\label{eq:eq22}
\end{equation}  
Given three elements $z_r$, $z_s$, and $z_t$ of the basis with their respective exponent vectors $\alpha_r$, $\alpha_s$, and $\alpha_t$, we define the ring multiplication rule as $z_t:=z_r\circ z_s$, where $t=\mathcal{J}(\alpha_r+\alpha_s)$, $r=\mathcal{J}(\alpha_r)$, and $s=\mathcal{J}(\alpha_s)$. We may consider a more general multiplication rule with respect to an arbitrary basis, allowing the ring multiplication to yield a linear combination of several basis elements
\begin{equation}
C_{rs}^t z_t:=z_r\circ z_s.
\label{eq:eq23}
\end{equation}   
We refer to $C_{rs}^t$ as structure constants\footnote{Note that, $C_{rs}^t$ are similar, but not equivalent, to the structure constants of Lie algebras \cite{Gilmore2012}. In fact, viewing $\mathds{R}[x]$ as an algebra, the structure constants of the algebra are identically zero because polynomial rings are commutative. However, $C_{rs}^t$ in Eq.~(\ref{eq:eq23}) are not all zeros.} of the ring. We can now substitute Eqs.~(\ref{eq:eq18}), (\ref{eq:eq19}), and (\ref{eq:eq21}) in Eq.~(\ref{eq:eq17}) to derive a matrix-form representation of infeasibility certificates
\begin{equation}
\begin{aligned}
-1=&\mathrm{P}^{;rs}z_r\circ z_s+\left[\mathrm{P}^{j_1;rs}z_r\circ z_s\right]\circ\left[\mathrm{G}^{p_1}_{j_1}z_{p_1}\right]+\\
		&\left[\mathrm{P}^{j_1j_2;rs}z_r\circ z_s\right]\circ\left[\mathrm{G}^{p_1}_{j_1}z_{p_1}\right]\circ\left[\mathrm{G}^{p_2}_{j_2}z_{p_2}\right]+\cdots+\\
		&\left[\mathrm{P}^{j_1\cdots j_m;rs}z_r\circ z_s\right]\circ\left[\mathrm{G}^{p_1}_{j_1}z_{p_1}\right]\circ\cdots\circ\left[\mathrm{G}^{p_m}_{j_m}z_{p_m}\right]+\left[\mathrm{B}^{ic}z_c\right]\circ\left[\mathrm{H}^{a}_{i}z_{a}\right],
\end{aligned}
\nonumber
\end{equation} 
from which it follows 
\begin{equation}
-1=\left[\mathrm{P}^{rs}\mathrm{U}^{t}_{rs}+\mathrm{P}^{j_1;rs}\mathrm{U}^{t}_{j_1;rs}+\cdots+\mathrm{P}^{j_1\cdots j_m;rs}\mathrm{U}^{t}_{j_1\cdots j_m;rs}+\mathrm{B}^{ic}\mathrm{V}^{t}_{ic}\right]z_t,
\label{eq:eq24}
\end{equation}   
where
\begin{equation}
\mathrm{U}^{t}_{j_1\cdots j_k;rs}:=C_{rs}^{t_0}C_{t_0p_1}^{t_1}\cdots C_{t_{k-1}p_k}^{t}\mathrm{G}^{p_1}_{j_1}\cdots\mathrm{G}^{p_k}_{j_k},\quad k=0\cdot\cdot m,
\label{eq:eq25}
\end{equation}  
\begin{equation}
\mathrm{V}^{t}_{ic}:=\mathrm{H}_{i}^{a}C_{ca}^{t},
\label{eq:eq26}
\end{equation} 
in which $\mathrm{P},\mathrm{P}^{j_1},\cdots,\mathrm{P}^{j_1\cdots j_m}\succeq 0$ and $\mathrm{B}\in\mathds{R}^{\ell\times s_p(\rho)}$ are unknown matrices. We have, thus, transformed Eq.~(\ref{eq:eq17}) into a feasibility problem of SDPs \cite{Ramana1997}, where we seek positive semidefinite matrices $\mathrm{P},\mathrm{P}^{j_1},\cdots,\mathrm{P}^{j_1\cdots j_m}$ and a linear matrix $\mathrm{B}$ that satisfy the equalities in Eq.~(\ref{eq:eq24}). In writing Eq.~(\ref{eq:eq24}), we adopted the Einstein summation convention, where repeating a dummy index as a superscript and subscript implies summation over an appropriate range of the index. Moreover, the symbols before semicolons in superscripts and subscripts are matrix indices, and those after semicolons are matrix-component indices. The summation rule applies to operations on both matrices and matrix components. Finally, we formulate Eq.~(\ref{eq:eq17}) as the following SDP feasibility problem (see \citet[Section 5.8]{boyd2004convex})
\begin{equation}
\begin{aligned}
& \min\limits_{\mathrm{P}^j,\mathrm{B}}\;0\\
& \quad \mathrm{s.t.}\\
& \sum_{j\in\mathrm{GEN}}\left\langle \mathrm{P}^j, \mathrm{U}_j^t\right\rangle+\left\langle \mathrm{B}, \mathrm{V}^t\right\rangle+\delta^{1t}=0, \quad t=1\cdot\cdot s_p(\rho),\\
&\mathrm{P}^{j}\succeq 0,\quad j=1\cdot\cdot|\mathrm{GEN}|, \quad \mathrm{B}\in\mathds{R}^{\ell\times s_p(\rho)},
\end{aligned}
\label{eq:eq27}
\end{equation} 
with $\mathrm{GEN}$ the index set of the generators in Eq.~(\ref{eq:eq18}), $\delta$ the Kronecker delta tensor, and $\left\langle A, B\right\rangle:=\mathrm{tr}(A^{\mathrm{T}}B)$ the standard matrix inner product. 

The structure constants defined by Eq.~(\ref{eq:eq23}) play an important role in the sparsity structure of the foregoing SDP relaxation problems. The function $\mathcal{J}(\alpha)$ and the structure constants $C_{rs}^t$ depend on the choice of the basis and monomial ordering for $\mathds{R}[x]$. Appropriate choices could significantly affect the sparsity and reduce the computation time \cite{Henrion2009, ahmadi2017improving}. However, there is no established procedure to determine the best choice for a given problem, although algorithms have been developed that can exploit the structure of SDP relaxations to efficiently solve these polynomial optimization problems \cite{gatermann2004symmetry,fawzi2015sparse,ahmadi2017improving}. Here, we derive explicit expressions for $\mathcal{J}(\alpha)$ and structure constants using a specific monomial ordering. Suppose that the standard monomial basis represented by $z$ are sorted according to the graded lexicographic order \cite[Chapter 2]{cox2006using} with respect to 
\begin{equation}
x_1<_{glex} x_2<_{glex}\cdots<_{glex} x_n.
\label{eq:eq28}
\end{equation} 
Then, the $\mathcal{J}$ function is given by
\begin{equation}
\mathcal{J}(\alpha)=\sum_{j=0}^{n}{n-j-\mathcal{A}_j(\alpha)-1\choose\mathcal{A}_j(\alpha)-1}, \quad \mathcal{A}_j(\alpha):=\sum_{i=j+1}^{n}\alpha_i, \quad j=0\cdot\cdot n,
\label{eq:eq29}
\end{equation} 
where we applied the convention $(-1)!:=+\infty$\footnote{This follows from the analogy between the factorial and $\Gamma$ function and the fact that $\Gamma(x)\to+\infty$ as $x\to0^+$.}. Using this index map, the multiplication table of $\mathds{R}[x]$ with respect to its standard monomial basis is readily constructed. The index map associated with monomial multiplication in $\mathds{R}[x]$ is written
\begin{equation}
\xi:\mathbb{N}\times\mathbb{N}\mapsto\mathbb{N}:(i,j)\mapsto k=\mathcal{J}(\mathcal{J}^{-1}(i)+\mathcal{J}^{-1}(j)),
\label{eq:eq30}
\end{equation}     
from which the structure constants can be constructed 
\begin{equation}
C_{ij}^k:=
\begin{cases}
1 \quad \mathrm{if} \quad k=\xi(i,j),\\
0 \quad \mathrm{if} \quad k\neq \xi(i,j).
\end{cases}
\label{eq:eq31}
\end{equation}

The SDP relaxations that we derived in this section are generally sparse because of sparse stoichiometry matrices and the form of thermodynamic inequalities. We leverage this sparsity structure to reduce the number of linear variables in Eq.~(\ref{eq:eq27}). Specifically, $\{\left\langle \mathrm{P}^j, \mathrm{U}_j^t\right\rangle=0,\;\forall j\}$ for most $t$, especially in the limit $t\to s_p(\rho)$. Let $t':=\{t>1:\mathrm{U}_j^t=0_{s_p(\rho)\times s_p(\rho)},\;\forall j \}$, $\mathcal{V}:=\mathrm{vec}(\mathrm{U}^{t'})$, and $\mathrm{b}:=\mathrm{vec}(\mathrm{B})$. Then, we have $\mathcal{V}\mathrm{b}=0$, so we can reduce the number of unknowns in $\mathrm{b}$ from $s_p(\rho)$ to $s_p(\rho)-|t'|$ using the decomposition $\mathrm{b}=\mathrm{b}_\mathrm{m}+N\omega$. Here, $\mathrm{b}_\mathrm{m}$ is the minimum-norm component of $\mathrm{b}$, $\omega\in\mathds{R}^{s_p(\rho)-|t'|}$, and $N$ is a matrix containing an orthogonal basis of $\mathrm{ker}(\mathcal{V})$. All these computations can be done offline before passing Eq.~(\ref{eq:eq27}) to SDP solvers. Of course, those $\mathrm{P}^j$, for which $\{\mathrm{U}_j^t=0_{s_p(\rho)\times s_p(\rho)},\;\forall t\}$, can also be eliminated from computations before solving Eq.~(\ref{eq:eq27}). 

\subsection{Global optimization}\label{sec:global}
Global optimization is an alternative approach to determine the feasibility of the CSS. We formulate a phase-I problem, in which feasibility is sought by minimizing a penalty function over a convex region obtained by relaxing the equality constraints. Recall, in the concentration space, the CSS is defined by intersecting an affine space corresponding to linear equality constraints and a nonconvex region associated with thermodynamic inequalities. It is more desirable to formulate this problem in the $Y$ space because the set of thermodynamically feasible concentrations is a convex polyhedron. Here, the penalty function measures the distance of any point in the space from the space defined by equality constraints, which is curved in the $Y$ space. Accordingly, nonconvexity is transferred from the constraints to the objective function, so feasibility is determined by minimizing a nonconvex function over a polyhedral set. This offers two major computational advantages: (i) Global solvers can efficiently partition the feasible set, construct and explore search trees, and identify the global solution through branch-and-bound or branch-and-reduce algorithms \cite{Tawarmalani2013}; (ii) Rounding errors in local solvers are alleviated when performing gradient-based search \cite{lasdon1978design} over distances that span several orders of magnitude in the concentration space. The phase-I problem corresponding to Eq.~(\ref{eq:eq27}) is
\begin{equation}
\begin{aligned}
& f^{\star}(\theta):=\min\limits_{y\in\mathds{R}^n}\|A\exp y-w-F\theta\|_2\\
& \quad \mathrm{s.t.}\\
& S^{\mathrm{T}}y \leq \kappa+\nu\ln\theta_1, \\
& y\leq 0,
\end{aligned}
\label{eq:eq32}
\end{equation}       
from which the infeasibility criterion 
\begin{equation}
\mathcal{Y}(\theta)=\emptyset \Leftrightarrow  f^{\star}(\theta)>0 
\label{eq:eq33}
\end{equation} 
follows. The search for the region of feasibility in the parameter space 
\begin{equation}
\Theta:=\left\{\theta\in \mathds{R}^{n_{\theta}}:\mathcal{Y}(\theta)\neq\emptyset\right\}
\label{eq:eq34}
\end{equation} 
is computationally expensive in general. However, because of the linear constraints in Eq.~(\ref{eq:eq10}), the search space can be reduced before solving Eq.~(\ref{eq:eq32}) by constructing a feasible set, only accounting for linear constraints
\begin{equation}
\mathcal{X}_{\mathrm{lin}}(\theta):=\left\{x\in [0,1]^{n}:Ax=w+F\theta, x\geq0\right\},
\label{eq:eq35}
\end{equation} 
\begin{equation}
\Theta_{\mathrm{lin}}:=\left\{\theta\in \mathds{R}^{n_{\theta}}:\mathcal{X}_{\mathrm{lin}}(\theta)\neq\emptyset\right\}.
\label{eq:eq36}
\end{equation}   
The set $\Theta_{\mathrm{lin}}\supseteq\Theta$ can be readily constructed using multiparametric linear programming algorithms \cite{Jones2007,Akbari2017} by solving the phase-I simplex problem \cite{Bertsimas1997}
\begin{equation}
\begin{aligned}
& f^{\star}_{\mathrm{lin}}(\theta):=\min\limits_{x,x'}\|x'\|_1\\
& \quad \mathrm{s.t.}\\
& Ax+x'=w+F\theta, \\
& x,x'\geq 0,
\end{aligned}
\label{eq:eq37}
\end{equation}  
providing a convex enclosure for $\Theta$.

\section{Characterization of the concentration solution space}\label{sec:characterize}
\begin{figure} [h]
\centering
\includegraphics[width=\linewidth]{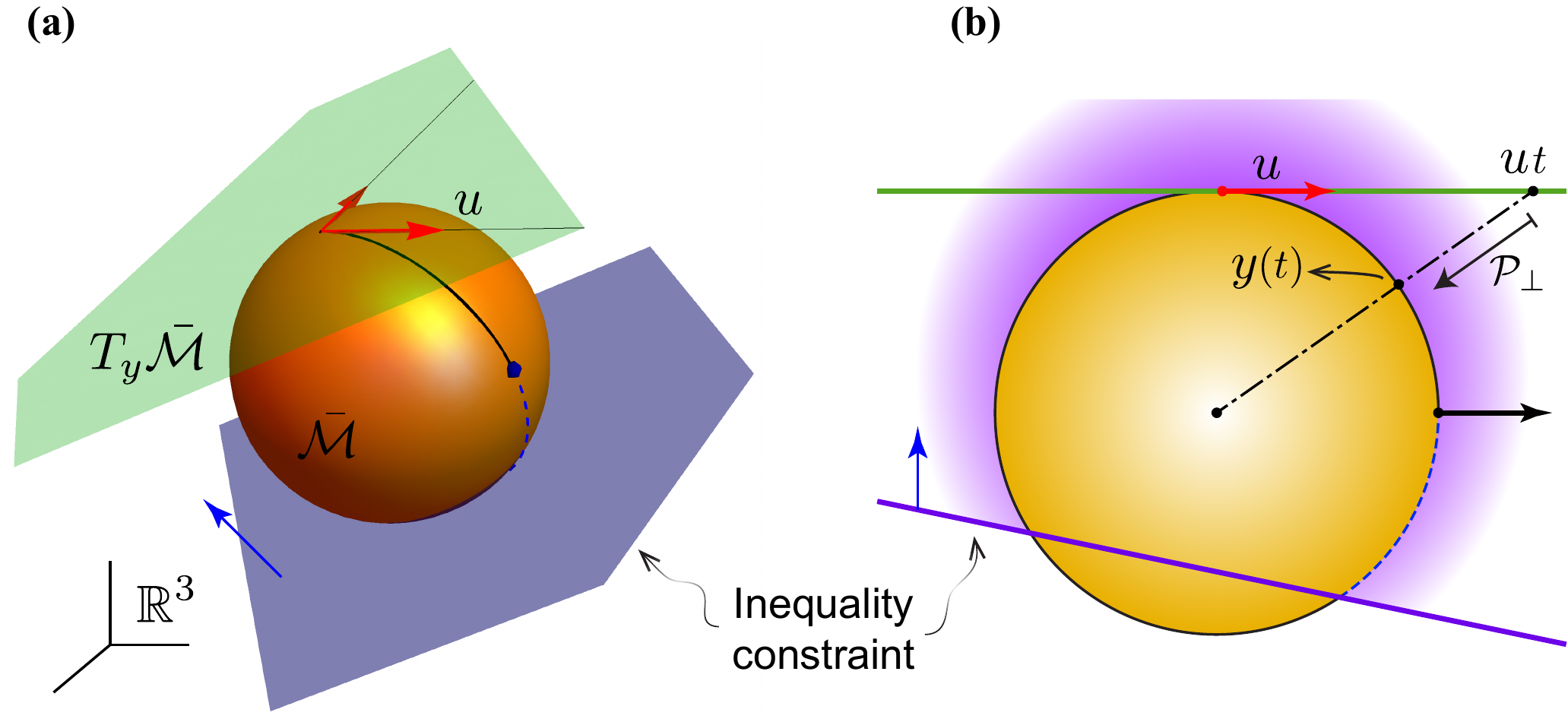}
\caption{Limitations of orthogonal projection in sampling the CSS when the curvature of the equality-constraint manifold $\bar{\mathcal{M}}$ is large. (a) Schematic representation of the sampling method. The CSS is explored by constructing curves in $\bar{\mathcal{M}}$ along random vectors $u$ generated in the tangent space $T_y\bar{\mathcal{M}}$. (b) A two-dimensional view of $\bar{\mathcal{M}}$, demonstrating part of the CSS that cannot be accessed through orthogonal projection (dashed line). Blue arrows indicate the feasibility direction of the inequality constraint. The geometry of $\bar{\mathcal{M}}$ is chosen here to illustrate the limitations of orthogonal projections, so it does not necessarily correspond to the exact geometry of manifolds arising from the equalities in Eq.~(\ref{eq:eq14}).}
\label{fig:fig2}
\end{figure}

Once the feasibility of $\mathcal{C}$ has been verified, it may be further analyzed to estimate characteristic parameters of biological relevance. In this section, we present two computationally tractable approaches to extract useful biological information from the CSS: (i) A determinist approach to compute global bounds on metabolite concentrations and reaction energies, and (ii) a stochastic approach to ascertain characteristic measures of concentrations and reaction energies. In the first, we solve 
\begin{equation}
y_{i,\mathrm{min}}:=\min_{y\in\mathcal{Y}} y_i,\quad i=1\cdot\cdot n,
\label{eq:eq38}
\end{equation} 
\begin{equation}
\Delta_rG'_{j,\mathrm{min}}:=\min_{y\in\mathcal{Y}} \Delta_rG'_j,\quad j=1\cdot\cdot m,
\label{eq:eq39}
\end{equation}  
to determine global lower bounds on metabolite concentrations and reaction energies. Upper bounds are calculated in a similar manner. In the second, we compute statistical characteristics of the CSS, such as expectations and standard deviations of concentrations and reaction energies. For a random variable $X:\mathcal{C}\mapsto\mathds{R}$ defined over the CSS, the expectation and standard deviation of $X$ are given by    
\begin{equation}
\mathbb{E}(X):=\int_{\mathcal{C}}X\mathrm{d}\mu,
\label{eq:eq40}
\end{equation} 
\begin{equation}
\sigma(X):=\sqrt{\mathbb{E}((X-\hat{X})^2)},
\label{eq:eq41}
\end{equation}   
where $\hat{X}:=\mathbb{E}(X)$, and $\mu$ is a probability measure of $X$. To simplify our analysis, we use the line measure to approximate these expectations and standard deviations. The line integrals are computed along random curves generated in the equality-constraint manifold (Fig.~\ref{fig:fig2}a).  

Here, we briefly outline the theoretical basis needed to construct random curves in the $Y$ space, along which all the equality constraints are satisfied. We denote the abstract $(n-\ell)$-dimensional manifold defined by the linear equalities in Eq.~(\ref{eq:eq14}) by $\mathcal{M}$ and its isometric embedding in $Y$ by $\bar{\mathcal{M}}$. We chose the embedding map $\iota:\mathcal{M}\hookrightarrow Y$ to simplify the construction and computation of the line measures. Moreover, let $\mathscr{A}:=A^{\mathrm{T}}$, $b:=(w+F\theta)^{\mathrm{T}}$, and 
\begin{equation}
\mathcal{F}^i(y):=\exp y^j\mathscr{A}^i_j-b^i.
\label{eq:eq42}
\end{equation} 
So far, we have not distinguished between subscripts and superscripts to identify the components of vectors and covectors. Vectors were viewed as $n$-tuples, the components of which were represented by subscripts (column vector in matrix form), irrespective of the space they were associated with. In this section, we denote the components of vectors and covector by superscripts and subscripts, which are represented by row and column vectors in matrix form. We apply the same convention to covariant and contravariant components of tensors. Given a full-rank $A$ matrix (\ie, $\mathrm{rank}(A)=\ell$ with $\ell<n$), we have
\begin{equation}
\bar{\mathcal{M}}:=\iota(\mathcal{M})=\{y\in\mathds{R}^n:\mathcal{F}(y)=0\}.
\label{eq:eq43}
\end{equation}  
Suppose $q\in\mathcal{M}$ and $y\in\bar{\mathcal{M}}$ such that $y=\iota(q)$. Then, the tangent space at $y$ is $T_y\bar{\mathcal{M}}=\mathrm{ker}(D\mathcal{F})$ \cite{ivancevic2007applied}. In matrix form, the derivative of $\mathcal{F}$ is given by $D\mathcal{F}=AE$ with $E:=\mathrm{diag}(\exp y)$. We derive explicit expressions for coordinate charts of $\mathcal{M}$ using null-space bases of $A$. Let $N\in\mathds{R}^{n\times(n-\ell)}$ denote a matrix, the columns of which form an orthonormal basis of $\mathrm{ker}(A)$. Noting that $\Phi=\iota\circ\varphi^{-1}$, the coordinate chart $\varphi$ at $y$ can be ascertained from $\Phi:\chi\mapsto \ln(x+\chi N^{\mathrm{T}})$ with $\chi\in\mathds{R}^{n-\ell}$. The components of $\chi$ should be viewed as the coordinates of $\varphi$ at $y$. Accordingly, an induced metric \cite{lee2006riemannian} for $\mathcal{M}$ is defined 
\begin{equation}
\begin{aligned}
g_{kl}&:=g_q\left(\frac{\partial}{\partial\chi^k},\frac{\partial}{\partial\chi^l}\right)=\bar{g}_y\left(\iota_{*}\left(\frac{\partial}{\partial\chi^k}\right),\iota_{*}\left(\frac{\partial}{\partial\chi^l}\right)\right)\\
&=\bar{g}_y\left(\mathcal{G}^i_k\frac{\partial}{\partial y^i},\mathcal{G}^j_l\frac{\partial}{\partial y^j}\right)=\mathcal{G}^i_k\mathcal{G}^j_l\delta_{ij},
\end{aligned}
\label{eq:eq44}
\end{equation} 
where $\bar{g}_{ij}:=\bar{g}_y(\partial/\partial y^i,\partial/\partial y^j)=\delta_{ij}$ is the standard Euclidean metric associated with $Y\equiv\mathds{R}^n$, and $\iota_{*}$ is the derivative of $\iota$. Therefore, $T_y\bar{\mathcal{M}}=\iota_{*}(T_q\mathcal{M})$, where $\mathcal{G}:=D\Phi^{\mathrm{T}}=N^{\mathrm{T}}E^{-1}$ in matrix form. 

To generate random curves in $\bar{\mathcal{M}}$, a feasible point of the CSS is required. It is preferably distanced far away from all the CSS boundaries in the $Y$ space to sample a more representative region. Interior points are generated by solving the following parametric problem
\begin{equation}
\begin{aligned}
&y_q(w):=\mathrm{arg}\max\limits_{r,y}\; r-w\|y\|_2,\\
& \quad \mathrm{s.t.}\\
&\exp y\mathscr{A}=b,\\
&y s_{(j)}+r\|s_{(j)}\|_2\leq \kappa_j+\nu_j\ln\theta_1,\quad j=\cdot\cdot m,
\end{aligned}
\label{eq:eq45}
\end{equation} 
where $s_{(j)}$ is the $j$th column of $S$. Note that, $y_q(0)$ is the Chebyshev center \cite[Section 8.5.1]{boyd2004convex} of the polyhedral set of thermodynamically feasible concentrations in the $Y$ space, which is restricted to $\bar{\mathcal{M}}$. When $w>0$, the objective function finds an interior point of the polyhedron, which is maximally distanced from all its boundaries with respect to the second norm, avoiding arbitrarily large negative $y$.

Several approaches can be adopted to generate random curves in $\bar{\mathcal{M}}$ from $y$. Here, we propose two methods. The first uses the orthogonal projection. First, a vector $u=u^k\partial/\partial\chi^k$ is generated in $T_q\mathcal{M}$ with respect to the basis $\partial/\partial\chi^k$, where the components $u^k\in[-1,1]$ are random numbers. The embedding of $u$ in $T_y\bar{\mathcal{M}}$ is ascertained using $\bar{u}=\iota_{*}(u)=\bar{u}^i\partial/\partial y^i$, in which $\bar{u}^i=u^k\mathcal{G}^i_k$. Next, a line $\bar{\mathcal{L}}$ passing through $y$ is defined in $T_y\bar{\mathcal{M}}$. It is constructed along the tangent vector $\bar{u}$ with the representation $\bar{\mathcal{L}}:=\{y'\in\mathds{R}^n:y'=y+\bar{u}t,t\in[0,t_{\mathrm{max}}]\}$. Finally, $\bar{\mathcal{L}}$ is orthogonally projected onto $\bar{\mathcal{M}}$ to construct a trajectory $\bar{\mathcal{T}}$ in $\bar{\mathcal{M}}$ by solving 
\begin{equation}
\begin{aligned}
&y'(t):=\mathrm{arg}\min\limits_{y''\in\mathds{R}^n}\; \frac{1}{2}\|y''-y-\bar{u}t\|_2,\\
& \quad \mathrm{s.t.}\\
&A\exp y''^{\mathrm{T}}=b^{\mathrm{T}}\\
\end{aligned}
\label{eq:eq46}
\end{equation} 
with 
\begin{equation}
\bar{\mathcal{T}}:=\{y'\in\mathds{R}^n: y' \ \mathrm{solves}\ \mathrm{Eq.~(\ref{eq:eq46})}, t\in[0,t_{\mathrm{max}}]\}.
\label{eq:eq47}
\end{equation} 
The Karush–Kuhn–Tucker (KKT) conditions \cite[Section 5.5.3]{boyd2004convex} for Eq.~(\ref{eq:eq46}) is 
\begin{equation}
\begin{cases}
y''+\lambda^{\mathrm{T}}AE=y+\bar{u}t,\\
A\exp y''^{\mathrm{T}}=b^{\mathrm{T}},\\
\end{cases}
\label{eq:eq48}
\end{equation} 
where $\lambda$ is the dual variable corresponding to the equality constraints in Eq.~(\ref{eq:eq46}). It is a column vector associated with the co-normal space $N_y^{*}\bar{\mathcal{M}}$. Furthermore, $b$ is associated with the normal space $N_y\bar{\mathcal{M}}$, so it is a row vector. All the other vectors are associated with $T_y\bar{\mathcal{M}}$ and are row vectors. The KKT system Eq.~(\ref{eq:eq48}) is differentiated with respect to $t$, yielding 
\begin{equation}
\begin{cases}
\dot{y}''+\dot{\lambda}^{\mathrm{T}}AE+\lambda^{\mathrm{T}}A\mathcal{E}\dot{y}''^{\mathrm{T}}=\bar{u},\\
AE\dot{y}''^{\mathrm{T}}=0\\
\end{cases}
\label{eq:eq49}
\end{equation}
with overdot denoting derivatives with respect to $t$, and $\mathcal{E}\in\mathds{R}^{n\times n\times n}$ is defined as 
\begin{equation}
\mathcal{E}_{k}^{ij}:=
\begin{cases}
\exp y''^i, \quad & i=j=k,\\
0, \quad & \mathrm{Otherwise}.
\end{cases}
\label{eq:eq50}
\end{equation}
We integrate Eq.~(\ref{eq:eq49}) with respect to $t$ subject to the initial conditions $y''(0)=y$ and $\lambda(0)=0$ until one of the thermodynamic constraints is violated. Note that, because $y\leq0$ is implicitly imposed by incorporating $\sum_{i=1}^{n}x_i=1$ as an equality constraint, it is automatically satisfied and needs not be checked directly along these trajectories.

The equality constraints in Eq.~(\ref{eq:eq10}) may yield a highly curved manifold in the $Y$ space, depending on the entries of $A$. Such large curvatures may render parts of the CSS inaccessible to orthogonal projection (see Fig.~\ref{fig:fig2}b). As a result, the integration of Eq.~(\ref{eq:eq49}) continues indefinitely without intersecting any thermodynamic constraints in finite $t$. Therefore, we propose a second method based on the Riemannian exponential map. Here, given the random vector $u$, we use the exponential map
\begin{equation}
\mathrm{Exp}_q:T_q\mathcal{M}\mapsto\mathcal{M}:u\mapsto\gamma(1)
\label{eq:eq51}
\end{equation}
to establish a correspondence\footnote{An isomorphism here.} between the line $\mathcal{L}\in T_q\mathcal{M}$ along the vector $u$ and geodesic $\gamma$. Unlike the orthogonal projection, the exponential map can access the entire $\mathcal{M}$ from $T_q\mathcal{M}$ for compact and connected manifolds, irrespective their curvature \cite{ivancevic2007applied}. Accordingly, trajectories are defined as $\mathcal{T}:=\{q'\in\mathcal{M}: q'=\mathrm{Exp}_q(ut), t\in[0,t_{\mathrm{max}}]\}$. Geodesics are curves of constant speed with respect to a connection on a manifold, that is they are solutions of $D_t\dot{\gamma}(t)=0$\footnote{$D_t$ is the covariant derivative of a curve parametrized with $t$.}. This allows to compute the components of $\gamma$ by solving a system ordinary differential equations. Let $\chi$ be a vector representing the coordinates of $\gamma$ with respect to the chart $\varphi$ we defined previously. Then, they are the solution of
\begin{equation}
\ddot{\chi}^k+\Gamma_{ij}^k\dot{\chi}^i\dot{\chi}^j=0,\quad k=1\cdot\cdot n-\ell,
\label{eq:eq52}
\end{equation}
subject to the initial conditions $\chi(0)=0$ and $\dot{\chi}(0)=u$. Here, $\Gamma_{ij}^k$ are the coefficients of an affine connection that $\mathcal{M}$ is equipped with. We require the connection to be compatible with the metric we derived in Eq.~(\ref{eq:eq44}), so we have $\Gamma_{ij}^k=g^{kr}(\partial_i g_{rj}+\partial_j g_{ri}-\partial_r g_{ij})/2$ \cite{lee2006riemannian}, where $\partial_s:=\partial/\partial \chi^s$. Although geodesics do not have the limitations of orthogonal projection, computing the trajectories $\mathcal{T}_i$ by solving Eq.~(\ref{eq:eq52}) is more expensive than by Eq.~(\ref{eq:eq49}). First, constructing a geodesic involves solving a system of coupled second-order differential equations. Second, integrating Eq.~(\ref{eq:eq52}) requires updating $(n-\ell)[(n-\ell)^2+3(n-\ell)+2]/6\sim O(n^3)$ unique, nonzero entries of $\Gamma_{ij}^k$\footnote{Note that, $\Gamma_{ij}^k$ is symmetric with respect to $i$ and $j$. It is also a dense matrix in general, so no sparsity structure is assumed at the outset.} at every time step compared to $2n\sim O(n)$ nonzero components of $E$ and $\mathcal{E}$ in Eq.~(\ref{eq:eq49}). A practical strategy is to perform initial numerical experiments with the orthogonal-projection method to estimate the fraction of trajectories, along which no thermodynamic constraint is violated when integrating Eq.~(\ref{eq:eq44}) for a prescribed $t_{\mathrm{max}}$. If this fraction does not meet a desired threshold, it implies that the CSS is located on a high-curvature region of $\bar{\mathcal{M}}$. Therefore, using a more expensive method, such a geodesic computations, may be necessary.             
 
Having constructed the random trajectories $\{\bar{\mathcal{T}}_i\}_{1}^{n_{\mathrm{traj}}}$, the expectation of a random variable $X$, defined by Eq.~(\ref{eq:eq40}), is approximated 
\begin{equation}
\mathbb{E}(X)\approx\frac{\sum\limits_{i=1}^{n_{\mathrm{traj}}}\int_{\bar{\mathcal{T}}_i}X\mathrm{d}l}{\sum\limits_{i=1}^{n_{\mathrm{traj}}}\int_{\bar{\mathcal{T}}_i}\mathrm{d}l},
\label{eq:eq53}
\end{equation}
where $n_{\mathrm{traj}}$ is the total number of randomly generated trajectories, and the line measure is given by $\mathrm{d}l=\sqrt{\left\langle \dot{x}',\dot{x}'\right\rangle}\mathrm{d}t=\sqrt{\left\langle \dot{y}'E,\dot{y}'E\right\rangle}\mathrm{d}t=\sqrt{\dot{y}'E^2\dot{y}'^{\mathrm{T}}}\mathrm{d}t$. The standard deviation of $X$ is approximated in a similar manner. 

\section{Numerical examples}\label{sec:example}
We present two case studies, demonstrating the application of the techniques we introduced in previous sections to characterize the CSS of biochemical reaction networks. In the first, we consider a small-scale toy model, for which polynomial-optimization, global-optimization, and random-sampling methods are all computationally tractable. We compare the computational performance of these approaches in determining the feasibility of the CSS. In the second, we study the CSS of glycolysis using global-optimization methods, illustrating their applicability to larger-scale pathways of biological relevance.  

\subsection{Toy model}\label{sec:toy}
\begin{figure} 
\centering
\includegraphics[width=\linewidth]{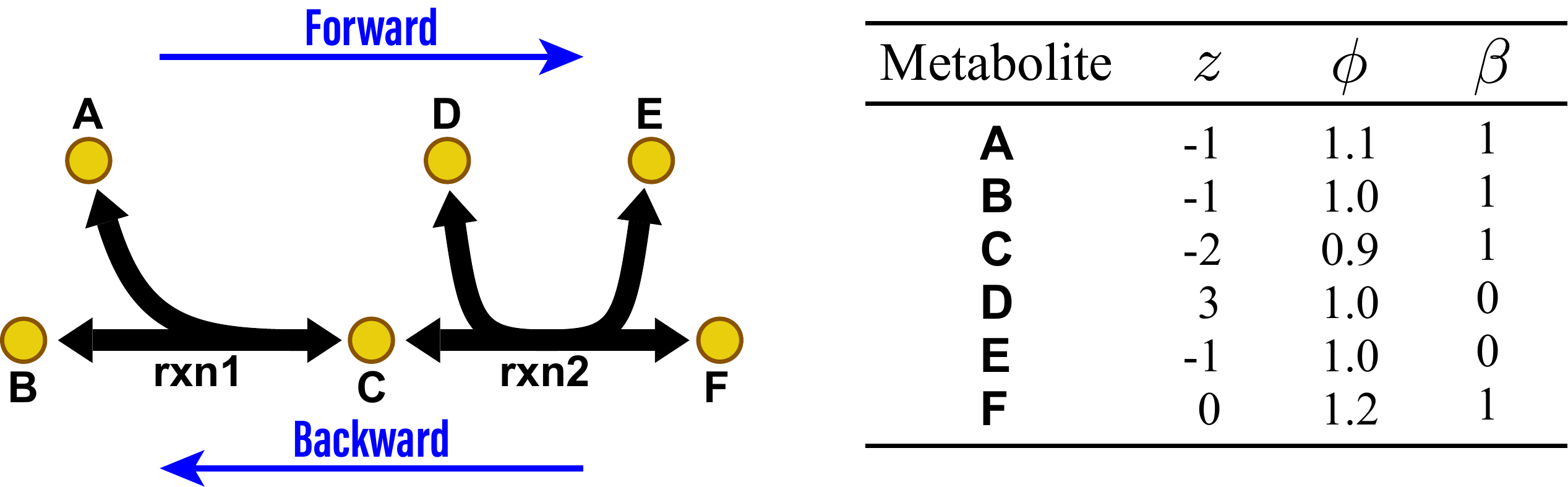}
\caption{Schematic representation of a toy model with two reactions and six metabolites. The apparent equilibrium constants for the two reactions are $K'_1=2.2\times10^4$ and $K'_2=1.22$. All stoichiometric coefficients are 1. The feasibility of the CSS is examined in the forward and backward direction.}
\label{fig:fig3}
\end{figure}

\begin{figure} [t]
\centering
\includegraphics[width=\linewidth]{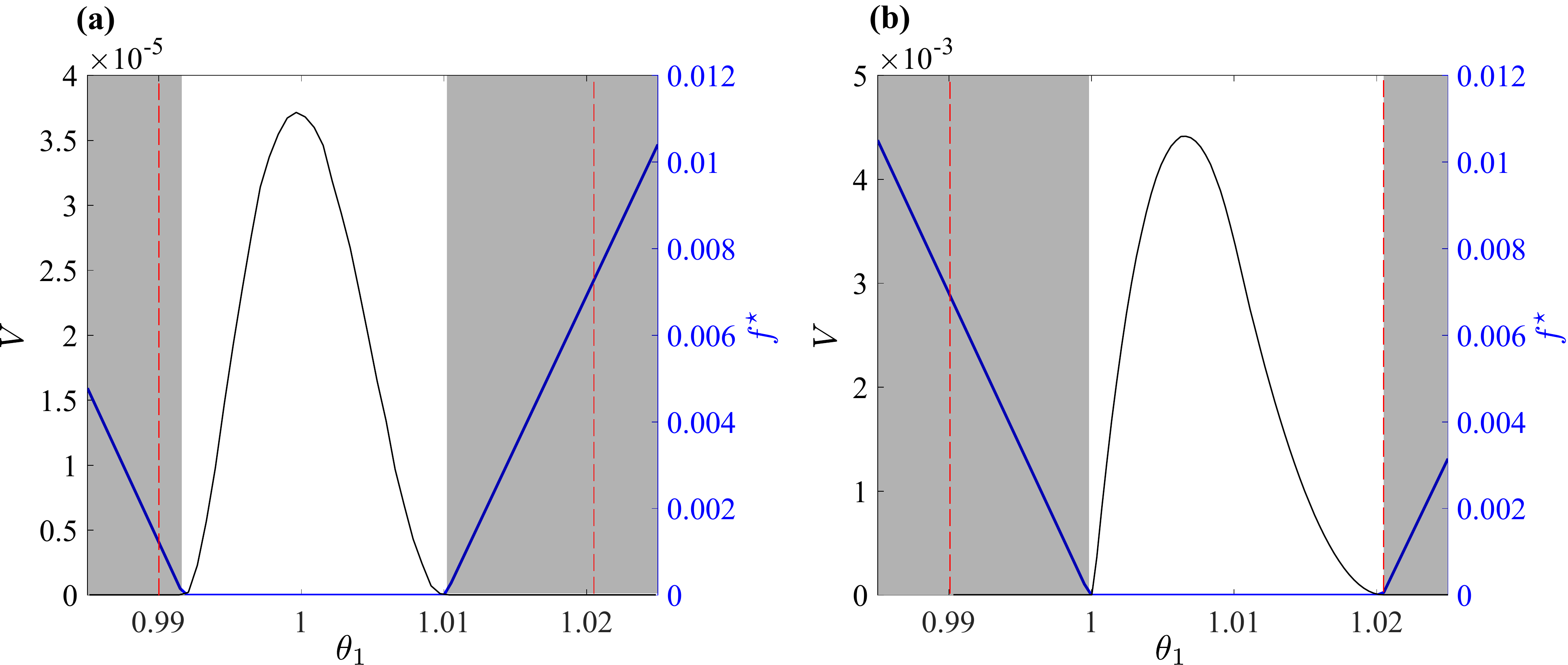}
\caption{Comparison of polynomial-optimization (shaded areas), global-optimization (blue lines), and random-sampling (black lines) methods in determining the emptiness of the CSS in the (a) forward and (b) backward direction of the toy model shown in Fig.~\ref{fig:fig3}. Parametric feasibility is determined with respect to $\theta_1$ along the line $\theta_2=0.1\theta_1$. Vertical dashed lines indicate a range of $\theta_1$, in which $\mathcal{X}(\theta)\neq\emptyset$ only when linear constraints are included. In the random-sampling method, the volume $V$ of the CSS in the affine space defined by the equality constraints in Eq.~(\ref{eq:eq10}) is used as a feasibility measure. At each grid point, the volume is estimated by generating $10^8$ random points with a uniform distribution in a bounding box. In the polynomial-optimization method, shaded areas indicate ranges of $\theta_1$, for which $\mathcal{X}(\theta)=\emptyset$ according to Eq.~(\ref{eq:eq27}).}
\label{fig:fig4}
\end{figure}

Consider the system of two hypothetical reactions with six metabolite shown in Fig.~\ref{fig:fig3}. We determine the parametric feasibility of this model with respect to $\theta_1$ at fixed $\mathcal{B}$ using polynomial-optimization, global-optimization, and random-sampling methods to compare their computational performance. The goal is to ascertain a range of $\theta_1$, in which $\mathcal{X}(\theta)\neq\emptyset$ along the line $\theta_2=(\mathcal{B}/C_s)\theta_1$. The feasible parameter region $\Theta$ computed by these methods for the forward and backward direction of the model are compared in Fig.~\ref{fig:fig4}. We determined the linear feasibility region $\Theta_{\mathrm{lin}}$ in the parameter space by solving Eq.~(\ref{eq:eq37}) to establish a search region for parameter-sweep computations. The search region was then discretized into 80 intervals, along which $\mathcal{X}(\theta)$ was evaluated. All three methods could estimate $\Theta$ with a reasonable accuracy for this model. Note that, the feasibility region estimated by the polynomial-optimization method must be regarded as an outer approximation of $\Theta$ because this method only provides infeasibility certificates. Nevertheless, the outer approximation tightly bounds $\Theta$ in both reaction directions using only two optimization levels ($d\leq2$). 

From a biological standpoint, this simple model reveals a crucial feature of biochemical reaction networks, namely the ability to reverse the direction of reactions without affecting the functional state of the cell. This manifests in the overlapping feasible regions of parameters in the forward and backward direction (Fig.~\ref{fig:fig4}). Observe that, in the range $1\lesssim\theta_1\lesssim 1.01$, $\mathcal{X}(\theta)\neq\emptyset$ in both directions. This indicates an optimal range for $C_{t,c}$, in which the direction of the reaction system in Fig.~\ref{fig:fig3} can be reversed without violating any biophysical constraint. Glycolysis is an important pathway, the function of which hings on this property. Several reactions in this pathway operate near their equilibrium, enabling the cell to rapidly reverse their flux direction in response to changes in carbon sources and nutrients \cite{Bennett2009,Park2019} with minimal disturbance to global parameters controlling cellular operations. The total metabolite concentration $C_{t,c}$ (related to ionic strength) and buffer capacity $\mathcal{B}$ are the global parameters of our model, which play a significant role in important cellular processes, such as enzyme activities \cite{liu2017ionic} and $\mathrm{pH}$ homeostasis \cite{slonczewski2009cytoplasmic}.    

\begin{figure} 
\centering
\includegraphics[width=\linewidth]{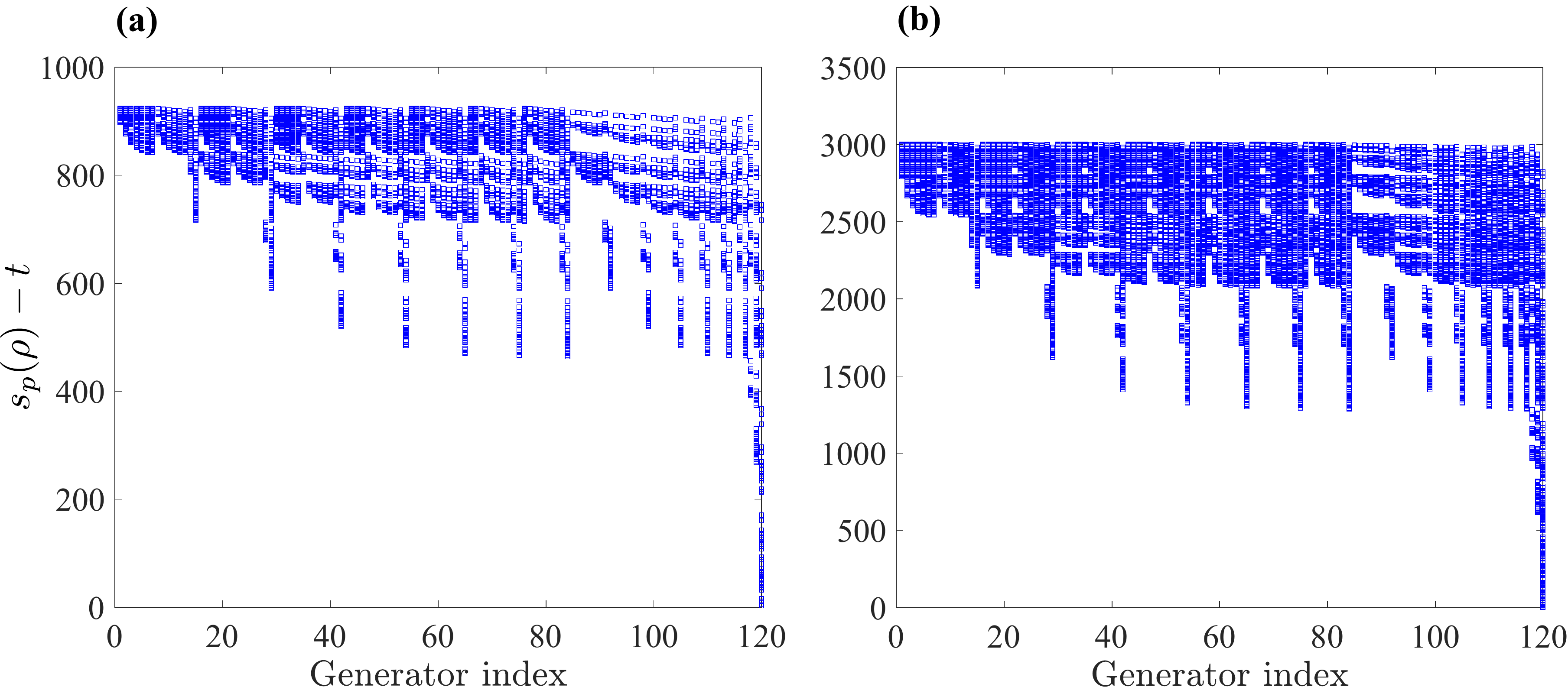}
\caption{Sparsity structure of the coefficient matrices $\mathrm{U}^t_j$ associated with the semidefinite matrices $\mathrm{P}^j$ in Eq.~(\ref{eq:eq27}) in (a) first- ($d=1$) and (b) second- ($d=2$) level optimization for the toy model in Fig.~\ref{fig:fig3}. Generator indices only include the first three terms in Eq.~(\ref{eq:eq18}). Symbols indicate coefficient matrices $\mathrm{U}^t_j$ with at least one nonzero entry.}
\label{fig:fig5}
\end{figure}

As previously stated, the sparsity of stoichiometry matrices in metabolic networks can be leveraged to lower the computation time of Eq.~(\ref{eq:eq27}). Figure~\ref{fig:fig5} shows the sparsity structure of the coefficient matrices $\mathrm{U}^t_j$ in the first and second level SDPs. Here, all the coefficient matrices $\mathrm{U}^t_j$ in 48.8\% (29.8\%) of the equality constraints of Eq.~(\ref{eq:eq27}) in the first (second) level SDP are all-zero matrices. As a result, these equalities only place constraint on the matrix $B$. By imposing these constraints on $B$ at the outset according to the procedure we outlined in Section~\ref{sec:pop}, they can be eliminated from Eq.~(\ref{eq:eq27}) before passing it to SDP solvers. This also reduces the number of linear variables in $B$.   

We compared the computational performance of polynomial-optimization (with and without sparse formulation), random-sampling, and global-optimization methods. The results are summarized in Table~\ref{tbl:tbl1}.  We formulated the global-optimization problem Eq.~(\ref{eq:eq32}) in the General Algebraic Modeling System \cite{Bussieck2004} and solved using the global solver BARON \cite{Tawarmalani2013} equipped with the local nonlinear programming solver CONOPT and the linear programming solver CPLEX. Random sampling and volume computations were performed in MATLAB R2019b. The SDP relaxations Eq.~(\ref{eq:eq27}) were formulated in MATLAB R2019b and solved using the SDP solver SeDuMi, which can handle sparse and dense inputs \cite{sturm1999using}. All computations were performed on a 64-bit Windows machine with a 3.1 GHz Intel Xeon CPU. 

Among these techniques, global optimization required the least amount of time to identify $\Theta$ in the search region. Although the sparse formulation of SDP relaxations resulted in 2--3 fold improvement in computational speed, the application scope of polynomial optimization is limited to small-scale problems. This is because the number and size of the unknown semidefinite matrices $\mathrm{P}^j$ in Eq.~(\ref{eq:eq27}) grow combinatorially with the number of inequalities and dimensions in Eq.~(\ref{eq:eq10}). These are directly related to the number of reactions and metabolites in the network. Therefore, given the current state of SDP solvers, application of polynomial optimization to genome-scale models is not encouraging. Random sampling also suffers from similarly poor scaling properties, although it performed faster than polynomial optimization for our toy model. Hence, it is not a promising candidate for larger-scale metabolic networks of practical interest, especially because the CSS is generally nonconvex and irregularly shaped in higher dimensions.       

\begin{table}
\centering
		\caption{Comparison of computation times between polynomial-optimization (PO), sparse polynomial-optimization (PO-sparse), random-sampling (RS), and global optimization (GO) methods to generate the parametric feasibility results in Fig.~\ref{fig:fig4} corresponding to the toy model shown in Fig.~\ref{fig:fig3}. Computation times are reported in seconds. In all cases, the search region in the parameter space was discretized into 80 intervals, and parameter-sweep computations were performed along the grid points.}
		\begin{tabular} {l c c c c}
		\hline 
		Direction & PO & PO-sparse & RS & GO \\
		\hline 
		Forward & 3248.4 & 1459.7 & 838.7 & 346.9\\
		Backward & 7644.8 & 2320.5 & 823.7 & 376.4\\
		\hline
		\end{tabular}		
		\label{tbl:tbl1}
\end{table} 

\subsection{Glycolysis}\label{sec:glycolysis}
\begin{figure} 
\centering
\includegraphics[width=\linewidth]{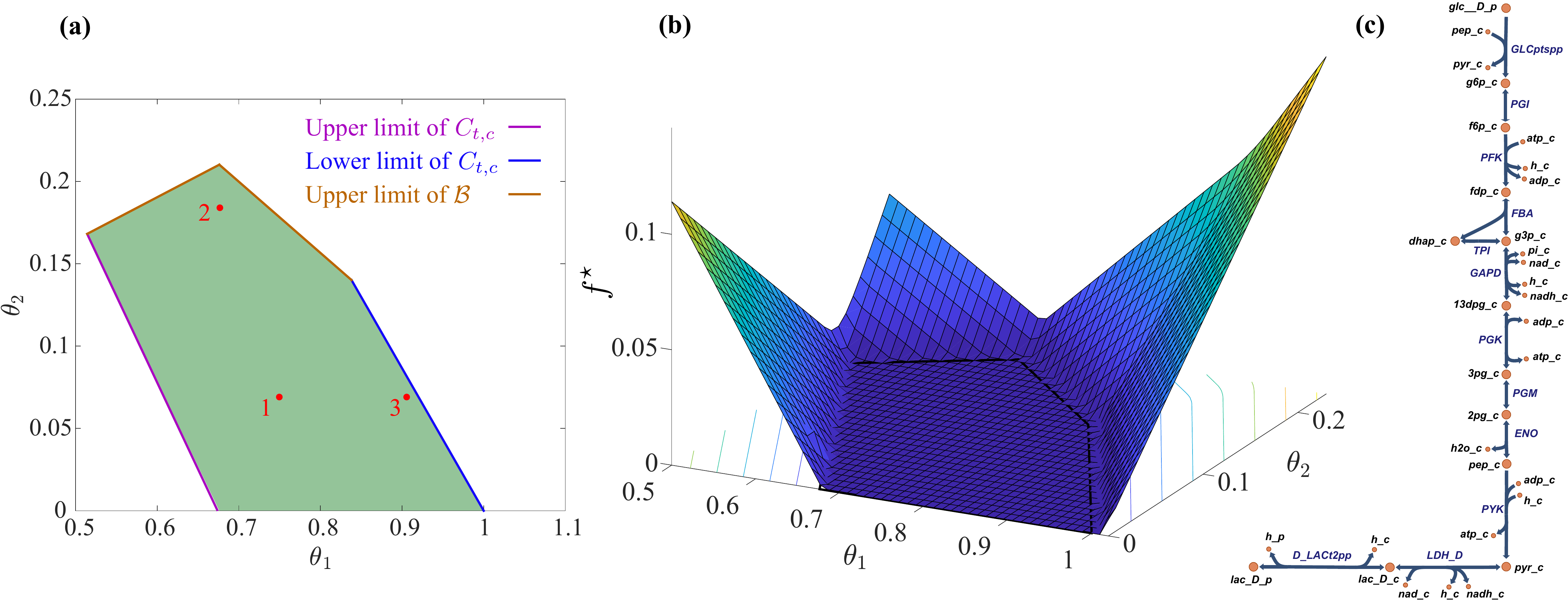}
\caption{Characterization of the CSS for glycolysis using global-optimization methods. (a) The linear feasibility region $\Theta_{\mathrm{lin}}$. Red circles indicate the points, at which the CSS is characterized in Fig.~\ref{fig:fig7}. (b) The solution of Eq.~(\ref{eq:eq10}), showing part of the parameter space, in which $\mathcal{X}(\theta)\neq\emptyset$. (c) Schematic representation of the glycolytic pathway, for which the CSS is analyzed.}
\label{fig:fig6}
\end{figure}

We construct the CSS for glycolysis, which comprise 21 metabolites and 12 reactions (Fig.~\ref{fig:fig6}c), and characterize it using the trajectory-tracing method outline in Section~\ref{sec:characterize}. For this model, both polynomial-optimization and random-sampling methods failed to identify the feasibility region $\Theta$ within prescribed time limits. Thus, we only present results for the global-optimization method in this section. 

We first constructed the linear feasibility region $\Theta_{\mathrm{lin}}$ for glycolysis (Fig.~\ref{fig:fig6}a) by solving Eq.~(\ref{eq:eq37}) using parametric linear-programming algorithms and identifying the critical regions, in which $f^{\star}(\theta)=0$ \cite{Akbari2017}. Because all the linear constraints are also linearly parametrized with respect to $\theta$, $\Theta_{\mathrm{lin}}$ is a polyhedron. This linearity results from the idealized model described in Section~\ref{sec:model}, where metabolite buffer capacities and charges are treated as constant. However, these coefficients can be functions of the parameters $\theta$ in more realistic models of biophysical constraints \cite{Akbaria2020}, resulting in a feasibility region with curved boundaries. We then solved Eq.~(\ref{eq:eq32}) to construct $\Theta$ within $\Theta_{\mathrm{lin}}$ and found that these regions coincide for glycolysis (Fig.~\ref{fig:fig6}b). This implies that the CSS is restricted by charge-related constraints as $\theta$ approaches the boundaries of $\Theta_{\mathrm{lin}}$ in Fig.~\ref{fig:fig6}a. Note that, $\Theta$ is a bounded region from above and below with respect to $\mathcal{B}$ and $C_{t,c}$, demonstrating why these global parameters cannot assume arbitrary values in biological systems that operate within the bounds of fundamental and evolutionary biophysical constraints.   

\begin{figure} 
\centering
\includegraphics[width=\linewidth]{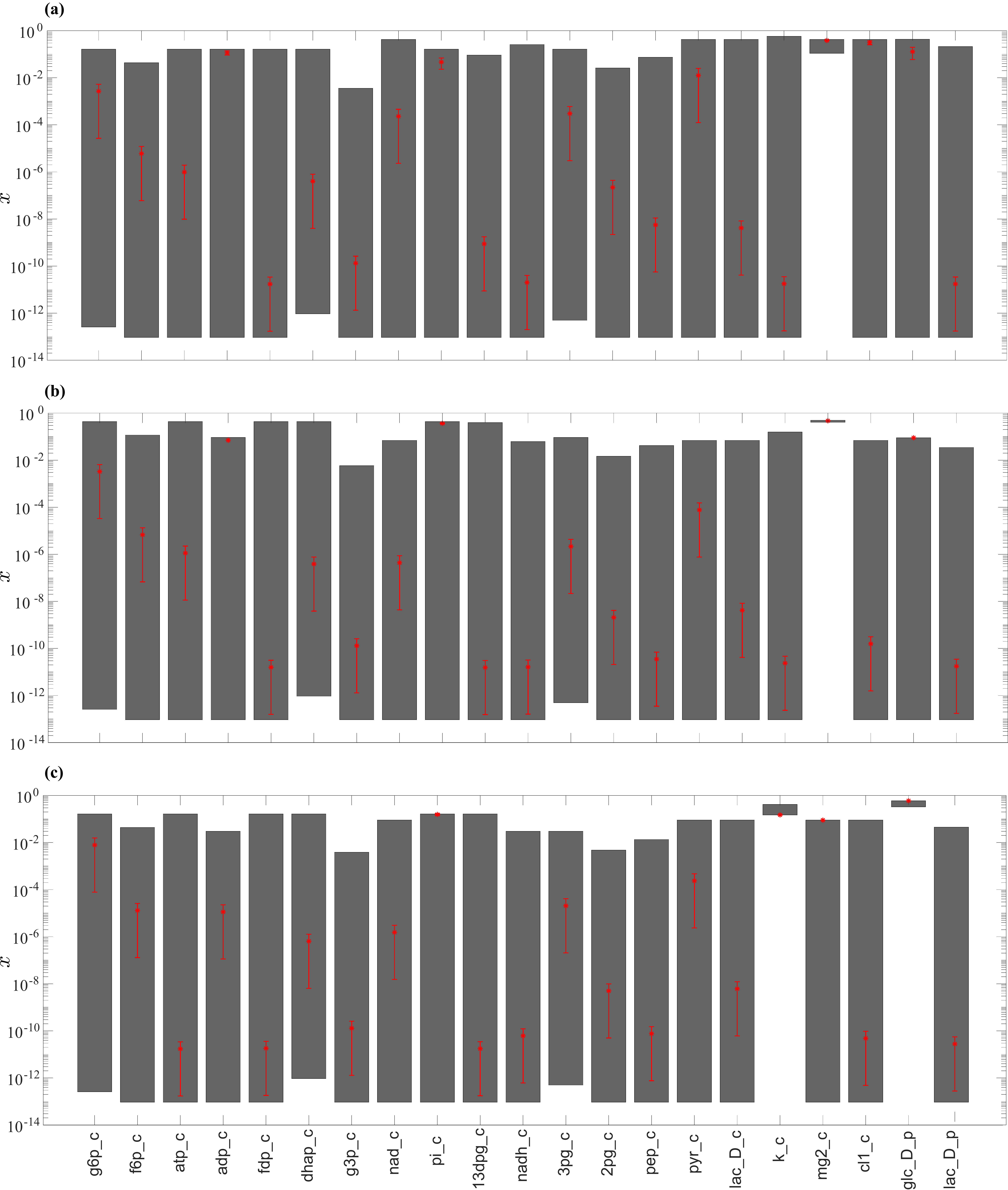}
\caption{Feasible concentration ranges of metabolites involved in glycolysis, which are evaluated at (a) $(\theta_1,\theta_2)=(0.750,0.069)$, (b) $(\theta_1,\theta_2)=(0.677,0.184)$, and (c) $(\theta_1,\theta_2)=(0.906,0.069)$, corresponding to the points 1,2, and 3 in Fig.~\ref{fig:fig6}. Asterisks and error bars indicate the expectation and standard deviation of concentrations estimated by constructing 1000 random trajectories in the CSS.}
\label{fig:fig7}
\end{figure}     

Furthermore, we characterized the CSS at three points inside $\Theta$ (red circles in Fig.~\ref{fig:fig6}a) by computing the expectation and standard deviations of metabolite concentrations (red asterisks and error bars in Fig.~\ref{fig:fig7}). These were determined by generating 1000 random curves inside the CSS at each parameter point through orthogonal-projection and geodesic computations described in Section~\ref{sec:characterize}. We performed several numerical experiments to estimate the computational performance of these methods. Although the difference between the expectations and standard deviations furnished by these techniques was not appreciable (less than 1\%), computation times for constructing geodesics were significantly larger than for orthogonal projections (more than 10 times). Moreover, all the trajectories constructed in these experiments using the orthogonal-projection method terminated upon intersecting a thermodynamic constraint, indicating that the CSS lies in a low-curvature part of $\bar{\mathcal{M}}$. 

We also computed global feasible concentration ranges for the metabolites participating in glycolysis (gray bars in Fig.~\ref{fig:fig7}). These provide useful information about how dominant constraints vary with respect to parameters, especially along the boundaries of the feasibility region $\partial\Theta$. Here, the biophysical constraints shrink the feasible concentration range of some the metabolites, resulting in $\mathcal{X}(\theta)=\emptyset$ on $\partial\Theta$. As previously stated, in the case of glycolysis, charge-related constraints (\eg, osmotic and charge balance) are dominant on $\partial\Theta$ because $\Theta_{\mathrm{lin}}\equiv\Theta$. Consider Point 3 in Fig.~\ref{fig:fig6}a, for example, which is close to the lower limit of $C_{t,c}$ on $\partial\Theta$. In this limit, the concentrations of K$^+$ (k\_c) and glucose (glc\_D\_p) are restricted to high values within a narrow range (Fig.~\ref{fig:fig7}c) since these have the largest $\phi_i$ among the metabolites in the system. Consequently, they can contribute the most to the cytoplasmic osmotic pressure to maintain a constant $\Delta\Pi$ at low $C_{t,c}$. In contrast, at Point 2 in Fig.~\ref{fig:fig6}a, which is close to the upper limit of $\mathcal{B}$, the concentration of Mg$^{2+}$ (mg2\_c) is driven to high levels (Fig.~\ref{fig:fig7}b). This is because only negatively charged metabolites in the model can contribute to the overall buffer capacity. As a result, the concentration of all the negatively charge  metabolites in this limit is high. On the other hand, the buffer intensity of positively charged metabolites, such as Mg$^{2+}$ and K$^+$  is zero. However, Mg$^{2+}$ is a more effective candidate to counterbalance the resulting negatively charged system than K$^+$  because it is a bivalent ion. Therefore, only high concentrations of Mg$^{2+}$ can maintain electroneutrality at fixed $C_{t,c}$.        

\section{Concluding remarks}\label{sec:conclusion}
We studied the concentration solution space of a constraint-based model of intracellular concentrations. It contains all possible metabolite concentrations satisfying the biophysical constraints imposed on biological networks. Understanding how the characteristics of an organism and its environment shape this set is key to quantitative models of cellular functions.   

We examined three computational approaches to determine the conditions, under which the concentration solution space is non-empty, namely polynomial-optimization, random-sampling, and global-optimization methods. In the polynomial optimization formulation, we provided a strategy to exploit the sparsity structure of the problem to reduce the computation time of infeasibility-certificate verification. In the global-optimization formulation, we introduced a technique to leverage the linear and logarithmically linear forms of the constraints to improve the computational efficiency of parametric feasibility studies. We then compared the scaling properties of these techniques by applying them in two case studies: (i) A toy model with 2 reactions and 6 metabolites, and (ii) a glycolytic pathway comprising 12 reactions and 21 metabolites. The global optimization formulation exhibited better scaling properties than the other two, and it was the only approach, with which the concentration solution space could be characterized in both case studies. Moreover, we recently demonstrated the computational viability of global optimization for characterizing intracellular concentrations in the core metabolism of \emph{Escherichia coli} \cite{Akbaria2020}, so it is currently the most promising approach for handling large-scale metabolic networks of practical interest. 

Polynomial-optimization methods are theoretically sound, but computationally challenging. Their application in systems biology is also uncommon. However, unlike global-optimization techniques, their performance does not depend on the convexity or connectedness of the feasible set. Infeasibility certificates can be computed to analyze these feasible sets, irrespective of their convexity, by solving semidefinite-programming relaxations. These are convex optimization problems, for which several efficient algorithms are available \cite{majumdar2019survey}. Recent advances in memory management and economic sparse arithmetic for semidefinite programs are encouraging \cite{yurtsever2019scalable}. Sparsity-exploiting algorithms are particularly relevant to the constraint-based model we discussed in this paper, because its semidefinite programming relaxations are highly sparse. This is due to (i) sparse stoichiometry matrices and (ii) algebraic form of the thermodynamic constraints, always comprising no more than two monomials. We proposed a strategy to partly leverage this sparsity structure during preprocessing. However, significantly more improvements in computational speed is expected by adapting primal-dual algorithms to a specific sparsity structure of these relaxations, such as those discussed by \citet{yurtsever2019scalable} and \citet{majumdar2019survey}. In light of these recent developments, the applicability of polynomial optimization to pathways larger than glycolysis (\eg, the entire central-carbon metabolism) in the near future is promising.    

\section*{Credit author statement}
 Conceptualization (A.A., B.O.P.); Data curation (A.A.); Formal analysis (A.A.); Funding acquisition (B.O.P.); Investigation (A.A.); Methodology (A.A.); Project administration (B.O.P.); Resources (B.O.P.); Supervision (B.O.P.); Validation (A.A.); Visualization (A.A.); Roles Writing – original draft (A.A.); Writing – review \& editing (B.O.P.). 

\section*{Declaration of Competing Interest}
The authors declare that they have no known competing financial interests or personal relationships that could have appeared to influence the work reported in this paper.

\section*{Acknowledgements}
This work was funded by the Novo Nordisk Foundation Grant Number NNF10CC1016517.

\bibliographystyle{model6-num-names}
\bibliography{mybibfile}







\end{document}